\documentclass[letterpaper, 10 pt, conference]{ieeeconf}  

\IEEEoverridecommandlockouts                              

\usepackage{graphicx} 
\usepackage{enumerate}
\usepackage{amsmath} 
\usepackage{amssymb}  

\usepackage{amsthm}
\newtheorem*{remark}{Remark}
\newcommand{\E}{\mathbb{E}}

\title{\LARGE \bf
Bridging a gap in Kalman filtering output estimation
\\with correlated noises or direct feed-through 
\\from process noise into measurements
}

\author{Ameet S. Deshpande$^{1}$
\thanks{$^{1}$Ameet Deshpande is with GE Global Research, Niskayuna, NY, USA.
        {\tt\small ameet.deshpande@gmail.com}}%
}

\begin{document}

\maketitle
\thispagestyle{empty}
\pagestyle{empty}

\begin{abstract}
Traditional statements of the celebrated Kalman filter algorithm focus on the estimation of state, but not the output. For any outputs, measured or auxiliary, it is usually assumed that the posterior state estimates and known inputs are enough to generate the minimum variance output estimate, given by $y_{n|n} = C x_{n|n} + D u_{n}$. Same equation is implemented in most popular control design toolboxes. It will be shown that when measurement and process noises are correlated, or when the process noise directly feeds into measurements, this equation is no longer optimal, and a correcting term of $H w_{n|n}\doteq H \E(w_{n}|z_{n})$ is needed in above output estimation. This natural extension can allow designer to simplify noise modeling, reduce estimator order, improve robustness to unknown noise models as well as estimate unknown input, when expressed as an auxiliary output. This is directly applicable in motion-control applications which exhibits such feed-through, such as estimating disturbance thrust affecting the accelerometer measurements. Based on a proof of suboptimality \cite{KalmanCounterExample}, this correction has been accepted and implemented in Matlab 2016 \cite{MatlabKalmanFnc2016}.
\end{abstract}

\section{Introduction}
Consider a discrete-time system in the canonical form with time-update equations at step $n$ as below:
{\small
\begin{align} \label{eq:DiscreteDyn}
\begin{split}
x_{n+1} &= A  x_{n} + B  u + G  w \\
y &= C  x_{n} + D  u + H  w\\
z &= C_{m}  x_{n} + D_{m}  u + H_{m}  w + v
\end{split}
\end{align}
}
where $x$ is the state vector, $u$ is a known input vector, $w$ is unknown input or process noise, $z$ are measured outputs affected by measurement noise $v$ and $y$ is the vector of auxiliary outputs (which may well contain as a subset the measured outputs without the measurement noise ). The process and measurement noises are white with following correlations. $\mathtt{E}(ww^{T}) = Q$, $\mathtt{E}(vv^{T}) = R$ and $\mathtt{E}(wv^{T}) = N$.

\subsection{Problem of Output Estimation}
The problem we wish to solve is to compute minimum-variance estimate of the state $x_{n}$, output $y_{n}$ and step-ahead prediction $x_{n+1}$, given a history of measurements $z_{i}, i=n, n-1, n-2, \dots$. This can be denoted by the short-hand $x_{n|n}$, $y_{n|n}$ and $x_{{n+1}|n}$ respectively. The state estimation and prediction part of this problem is at the heart of linear real-time estimation of dynamic systems and was solved exactly by Kalman \cite{Kalman1960}, and is described in many texts  (e.g. see \cite{kwakernaak1972linear}, \cite{lewis1986optimal}, \cite{FranklinWorkmanPowell1997}, \cite{kailath2000}, \cite{simon2006}) with varying degree of generality of assumptions. The most general form, described for example in Kailath et al. \cite{kailath2000} has $G=I$, $H_{m}=0$ but allows arbitrary cross-correlation between $w$ and $v$ which makes it equivalent with \eqref{eq:DiscreteDyn} without loss of generality.

Yet there is a gap in all but one of above references. They do not explicitly describe the equations for optimal output estimation $y_{n|n} \doteq \E(y_{n}|z_{n})$. The only reference from those surveyed, which describes the output equation is Kwakernaak, Sivan \cite{kwakernaak1972linear} (see eq. 4-228, Thm. 4.7), and the output estimate equation is given by:
\begin{align}\label{eq:IncorrectOutputEq}
y_{n|n} = C x_{n|n} + D u_{n} 
\end{align}
This is proved for uncorrelated measurement and process noises, but a footnote states that the same can be proved for correlated noises, and that $\mathtt{var}(y - y_{n|n}) = H_{m}'Q H_{m} + R + H_{m}N + NH_{m}$, and innovations $y-y_{n|n}$ are white, when $y_{n|n}$ is computed as above and when the output vector is same as measurements  without the measurement noise, $y \doteq z-v$.

The steady-state or time-varying Kalman filter implementations in popular control design toolboxes such as Matlab 2015b \cite{MatlabKalmanFnc2015b}, Labview  \cite{LabviewKalmanFnc}, Mathematica \cite{MathematicaKalmanFnc}, Maple \cite{MapleKalmanFnc} and Octave-forge \cite{OctaveKalmanFnc} do allow us to pose the most general problem \eqref{eq:DiscreteDyn} with correlated noises or feed-through and also allow us to generate estimates for both states $x_{n|n}$ and outputs $y_{n|n}$. But as recently as 2015, all above implementations use the equation \eqref{eq:IncorrectOutputEq} for updating outputs. As a specific example, equation for {\emph{current}} form of Kalman filter per the documentaion in \cite{MatlabKalmanFnc2015b} is:
{\small
\begin{align}
\begin{bmatrix} x_{n|n}\\ y_{n|n}\\ x_{{n+1}|n}\end{bmatrix} &= 
\begin{bmatrix} 
	I - K_{g}C_{m} & -K_{g}D_{m}\\
     C - C K_{g} \,C_{m} & D - C K_{g} \,D_{m} \\
     A - M_{A,G} \,C_{m} & B - M_{A,G} \,D_{m}
\end{bmatrix} \cdot  \label{eq:XYXnextupdate:A}\\
& \hskip 6em 
\begin{bmatrix} x_{n|{n-1}}\\ u_{n}\end{bmatrix} +
\begin{bmatrix} K_{g}\\ C K_{g}\\ M_{A,G} \end{bmatrix} z_{n} \nonumber
\end{align}
}
where $M_{C} = C K_{g}$ and $M_{A,G} = A K_{g} + G K_{g2}$. $K_{g} = P_{n|{n-1}}C_{m}' (C_{m}P_{n|{n-1}}C_{m}' + \bar{R})^{-1}$, $K_{g2} = (Q H_{m}'+N) (C_{m}P_{n|{n-1}}C_{m}' + \bar{R})^{-1}$, $\bar{R} =  R + H_{m} Q H_{m}^{T} + H_{m} N + N^{T} H_{m}^{T}$ and $P_{n|n-1}$ is the solution of Riccati difference iteration or the Riccati equation. Despite the complex form, it is easy to see that eq. \eqref{eq:XYXnextupdate:A} matches eq. \eqref{eq:IncorrectOutputEq}. This paper proposes a change only to the output update equation for $y_{n|n}$.

As will be shown later in this paper, all claims from  \cite{kwakernaak1972linear} mentioned above hold except the optimality, and equation \eqref{eq:IncorrectOutputEq} fails to be optimal when there is a direct feed-through, which means either $H_{m}$ or $N$ is non-zero. In such case the correction needed is:
\begin{align} \label{eq:CorrecetdOutputEqA}
y_{n|n} = C x_{n|n} + D u_{n} + H \E(w_{n}|z_{n})
\end{align}

In section \ref{sec:Example}, we will prove the suboptimality of eq. \eqref{eq:CorrecetdOutputEqA} over eq. \eqref{eq:IncorrectOutputEq} through a counter-example.
The loss of optimality occurs because the conditional expectation of process noise $w_{n}$ knowing measurements $z_{n}$ is non-zero, as they are both correlated due to direct feed-through via $H_{m}$ or $N$, and must be accounted for to achive best estimation of output $y$. Another way to see this is in frequency domain. If \eqref{eq:DiscreteDyn} is a discretization of a continuous time system, then $x$ is bound to be smooth by physics, and so is the estimate $x_{n|n}$. But outputs $y$ is not smooth and would not have high-frequency roll-off, due to direct feed-through in $H$ or $H_{m}$. But using eq. \eqref{eq:IncorrectOutputEq} would constraint $y_{n|n}$ to be smooth even though $y$ is not. Thus estimate error is bound to be large at high frequencies. This can be avoided using the measurements $z_{n}$ in a better way, in the computation of $y_{n|n}$, as per eq. \eqref{eq:CorrecetdOutputEq} to follow.

\subsection{Problem of Unknown Input Estimation}
Note that after substituting $C=0$, $D=0$ and $H=I$ in \eqref{eq:DiscreteDyn}, we have $y=w$ and the problem of unknown input estimation becomes a specific case of output estimation.
Traditional methods of unknown input modelling include augmenting state-space to include noise states. Such methods work well for smoothly varying unknown inputs, but  due to the smoothness requirement on noise states evolution, they cannot model broad-band process noise which have significant energy at high frequencies or those which have a hard floor in frequency instead of a roll-off. Moreover, an estimator designed with a smooth noise model can lack robustness to unmodeled spikes in unknown inputs at high frequencies.
Especially for an application like a wind turbine, direct feed-through exists from unknown wind thrust to sensors like accelerometers. Due to pockets of localized wind gusts and tower dam effect, the thrust on blades can have a very broad-band spectrum extending to high-frequencies with peaks at multiples of blade passing frequency. 
The method to improve output estimation, proposed in this paper, has a fortunate side-effect of improving the unknown input estimation, and a direct way to model broad-band noise without expanding state-space. 

The paper is organized as follows. In sec. \ref{sec:DeriveDiscKalman}, we derive the discrete time minimum variance estimator or the Kalman filter from the first principles, with a specific goal of state as well as output estimation. Note that we only derive state update equation for the sake of completeness, and propose no change in them compared to prior art. The only change proposed is in output update equations. The comparison with earlier approaches for output estimation is done in sec. \ref{sec:ComparePrevMethod}. A simple numerical example which demonstrates the improved estimation is described in sec. \ref{sec:Example}.



\section{Derivation of discrete time Kalman filter}\label{sec:DeriveDiscKalman}
This section derives the formulae for a time-varying discrete-time Kalman filter from first principles. In doing so, we will leverage a key result from conditional probability of a bi-variate Gaussian distribution repeatedly.

\subsection{A result from conditional bivariate Gaussian distribution}\label{sec:bivarGaussDist}
Let us denote a normal distribution by its mean and variance, $\mathcal{N}(\mathtt{mean}, \mathtt{variance})$. Assume two normal random variables $x_{1} = \mathcal{N}(m_{1},P_{11})$ and $x_{2} = \mathcal{N}(m_{2},P_{21})$ with cross-covariance $\mathtt{cov}(x_{1}, x_{2})=\E((x_{1}-m_{1})(x_{2}-m_{2})^{T}) = P_{12}$. Then it is well-known \cite{CondBivarDistrib} that condition distribution of $x_{1}$ knowing $x_{2} = z_{2}$ is $x_{1|2} = \mathcal{N}(m_{1|2}, P_{1|2})$ where
{\small
\begin{align}
&m_{1|2} = \E(x_{1}|x_{2}=z_{2}) \nonumber\\
&\hskip1em = \E(x_{1}) + \mathtt{cov}(x_{1},x_{2}) \mathtt{var}(x_{2},x_{2})^{-1}(z_{2}-\E(x_{2})) \nonumber\\
&\hskip1em= m_{1} + P_{12} P_{22}^{-1} (z_{2} - m_{2})  \label{eq:condProb}\\
&P_{1|2} = \mathtt{var}(x_{1}|x_{2}=z_{2})\nonumber\\
&\hskip1em= \mathtt{var}(x_{1},x_{2}) - \mathtt{cov}(x_{1},x_{2}) \mathtt{var}(x_{2},x_{2})^{-1} \mathtt{cov}(x_{1},x_{2})' \nonumber\\
&\hskip1em= P_{11} - P_{12} P_{22}^{-1} P_{12}^{T}  \label{eq:condProbVar}
\end{align}
}
\subsection{State and Output update rules}\label{sec:updateEqs}
At the n'th step of filter update, we know prior state estimate $x_{n|{n-1}}$ (with variance $P_{n|{n-1}}$), measurements $z_{n}$ and known inputs $u_{n}$. Tabulating cross-covariances and variances will help us in upcoming development.
{\small
\begin{align*} 
&\mathtt{cov}(x_{n|{n-1}}, z_{n|{n-1}}) = P_{n|{n-1}} C_{m}^{T} \\
&\mathtt{cov}(y_{n|{n-1}}, z_{n|{n-1}}) = C P_{n|{n-1}} C_{m}^{T} + H Q H_{m}^{T} + H N\\
&\mathtt{cov}(x_{n+1|{n-1}}, z_{n|{n-1}}) = A P_{n|{n-1}} C_{m}^{T} + G Q H_{m}^{T} + G N\\
&\mathtt{var}(x_{n|{n-1}}) = P_{n|{n-1}}\\
&\mathtt{var}(z_{n|{n-1}}) = C_{m} P_{n|{n-1}} C_{m}^{T} + \bar{R}\\
&\mathtt{var}(x_{n+1|{n-1}}) = A P_{n|{n-1}} A^{T} + G Q G^{T} 
\end{align*}
}
where $\bar{R} =  R + H_{m} Q H_{m}^{T} + H_{m} N + N^{T} H_{m}^{T}$.
Using above expressions and equation \eqref{eq:condProb}, by variable substitution, posterior expectations and variances of various random variables can be readily computed. The posterior expectation of state is given by (measurement update for state) :
{\small
\begin{align}
&\E(x_{n|n}) = \E(x_{n|{n-1}} | z_{n|{n-1}} = z_{n}) \label{eq:Xupdate}\\
&= \E(x_{n|{n-1}}) + \mathtt{cov}(x_{n|{n-1}}, z_{n|{n-1}}) \cdot \nonumber\\
&\hskip 4em \mathtt{var}(z_{n|{n-1}})^{-1} (z_{n} - \E(z_{n|{n-1}})) \nonumber\\
&= x_{n|{n-1}} + P_{n|{n-1}} C_{m}^{T}  \cdot \nonumber\\
&\hskip 4em \left(C_{m} P_{n|{n-1}} C_{m}^{T} + \bar{R} \right)^{-1} \left(z_{n} - (C_{m} x_{n|{n-1}} +D_{m} u_{n}) \right) \nonumber
\end{align}
}
Similarly, the posterior expectation of output is given by (measurement update for outputs):
{\small
\begin{align}
&\E(y_{n|n}) = \E(y_{n|{n-1}} | z_{n|{n-1}} = z_{n}) \label{eq:Yupdate}\\
&= \E(y_{n|{n-1}}) + \mathtt{cov}(y_{n|{n-1}}, z_{n|{n-1}}) \cdot \nonumber\\
&\hskip 4em \mathtt{var}(z_{n|{n-1}})^{-1} (z_{n} - \E(z_{n|{n-1}})) \nonumber\\
&= C x_{n|{n-1}} + D u_{n} + \left(C P_{n|{n-1}} C_{m}^{T} + H Q H_{m}^{T} + H N \right) \cdot \nonumber\\
& \hskip 4em \left(C_{m} P_{n|{n-1}} C_{m}^{T} + \bar{R} \right)^{-1} \left(z_{n} - (C_{m} x_{n|{n-1}} +D_{m} u_{n}) \right)  \nonumber
\end{align}
}
and the posterior expectation of predicted state becomes (measurement update for predicted state):
{\small
\begin{align}
&\E(x_{{n+1}|n})  = \E(x_{n+1|{n-1}} | z_{n|{n-1}} = z_{n})  \label{eq:XupdateNext}\\
&= \E(x_{n+1|{n-1}}) + \mathtt{cov}(x_{n+1|{n-1}}, z_{n|{n-1}}) \cdot \nonumber\\
& \hskip 4em \mathtt{var}(z_{n|{n-1}})^{-1} (z_{n} - \E(z_{n|{n-1}})) \nonumber\\
&= A x_{n|{n-1}} + B u_{n} + \left(A P_{n|{n-1}} C_{m}^{T} + G Q H_{m}^{T} + G N \right) \cdot \nonumber\\
& \hskip 4em \left(C_{m} P_{n|{n-1}} C_{m}^{T} + \bar{R} \right)^{-1} \left(z_{n} - (C_{m} x_{n|{n-1}} +D_{m} u_{n}) \right) \nonumber
\end{align}
}
and the posterior variance for predicted state using eq. \eqref{eq:condProbVar} is given by the discrete Riccati differential equation below:
{\small
\begin{align}
&P_{{n+1}|n} = \mathtt{var}(x_{n+1|n-1}) - \mathtt{cov}(x_{n+1|n-1}, z_{n|n-1}) \cdot \label{eq:VarupdateXNext}\\
&\hskip 6em \mathtt{var}(z_{n|n-1})^{-1} \mathtt{cov}(x_{n+1|n-1}, z_{n|n-1})^{T} \nonumber\\
&=(A P_{n|{n-1}} A^{T} + G Q G^{T}) - \left(A P_{n|{n-1}} C_{m}^{T} + G Q H_{m}^{T} + G N \right) \cdot \nonumber\\
& \hskip 1em \left(C_{m} P_{n|{n-1}} C_{m}^{T} + \bar{R} \right)^{-1} 
\left(A P_{n|{n-1}} C_{m}^{T} + G Q H_{m}^{T} + G N \right)^{T} \nonumber
\end{align}
}
Similarly the posterior variance for predicted outputs is:
{\small
\begin{align}
&\mathtt{var}(y_{n|n}) = \mathtt{var}(y_{n|n-1}) - \mathtt{cov}(y_{n|n-1}, z_{n|n-1}) \cdot \label{eq:VarUpdateY}\\
&\hskip 6em \mathtt{var}(z_{n|n-1})^{-1} \mathtt{cov}(Y_{n|n-1}, z_{n|n-1})^{T} \nonumber\\
&= (C P_{n|{n-1}} C^{T} + H Q H^{T}) - \left(C P_{n|{n-1}} C_{m}^{T} + H Q H_{m}^{T} + H N \right) \cdot \nonumber\\
& \hskip 1em \left(C_{m} P_{n|{n-1}} C_{m}^{T} + \bar{R} \right)^{-1} 
\left(C P_{n|{n-1}} C_{m}^{T} + H Q H_{m}^{T} + H N \right)^{T} \nonumber
\end{align}
}

In summary, after each measurement $z_{n}$, $x_{n|{n-1}}$ and $P_{n|n+1}$ are propagated forward in time to generate $x_{{n+1}|n}$ and $P_{{n+1}|n}$ and recursively there on, through equations \eqref{eq:XupdateNext} and \eqref{eq:VarupdateXNext}. Estimates of outputs and states, $y_{n|n}$ and $x_{n|n}$ based on Kalman filter are also generated by equations \eqref{eq:Yupdate} and \eqref{eq:Xupdate}.

\begin{remark}
Note that if estimated outputs are same as measured, i.e. when $C=C_{m}$, $H=H_{m}$, $N=0$, eq. \eqref{eq:VarUpdateY} simplifies to 
\begin{align}
&\mathtt{var}(y_{n|n}) = (C_{m}P_{n|n-1}C_{m}^{T} +H_{m}Q H_{m}^{T}) \cdot  \label{eq:VarupdateYmeas}\\
& \hskip 3em (C_{m}P_{n|n-1}C_{m}^{T} +H_{m}Q H_{m}^{T}+R)^{-1} \cdot R \preceq R \nonumber
\end{align}
Thus the posterior output estimate is more accurate than the measurement, which matches the intuition.
\end{remark}

\section{Comparison with the prior art }\label{sec:ComparePrevMethod}
Equations \eqref{eq:Xupdate}, \eqref{eq:Yupdate} and \eqref{eq:XupdateNext} can be summarized in following matrix equation.
{\small
\begin{align}
\begin{bmatrix} x_{n|n}\\ y_{n|n}\\ x_{{n+1}|n}\end{bmatrix} &= 
\begin{bmatrix} 
	I - K_{g}C_{m} & -K_{g}D_{m}\\
     C - M_{C,H} \,C_{m} & D - M_{C,H} \,D_{m} \\
     A - M_{A,G} \,C_{m} & B - M_{A,G} \,D_{m}
\end{bmatrix} \cdot  \label{eq:XYXnextupdate}\\
& \hskip 6em 
\begin{bmatrix} x_{n|{n-1}}\\ u_{n}\end{bmatrix} +
\begin{bmatrix} K_{g}\\ M_{C,H}\\ M_{A,G} \end{bmatrix} z_{n} \nonumber
\end{align}
}
where $M_{C,H} = C K_{g} + H K_{g2}$ and $M_{A,G} = A K_{g} + G K_{g2}$. $K_{g} = P_{n|{n-1}}C_{m}' (C_{m}P_{n|{n-1}}C_{m}' + \bar{R})^{-1}$ and $K_{g2} = (Q H_{m}'+N) (C_{m}P_{n|{n-1}}C_{m}' + \bar{R})^{-1}$.

Using \eqref{eq:XYXnextupdate}, it can be shown that $y_{n|n}$ and $x_{n|n}$ are related in following way.
{\small
\begin{align} 
y_{n|n} &= C x_{n|n} + D u_{n} \label{eq:RelateYnnXnn}\\
&\hskip 2em  + H K_{g2} \cdot \left(z_{n} - ( C_{m} x_{n|{n-1}} + D_{m} u_{n} ) \right) \nonumber
\end{align}
}
Comparing eq. \eqref{eq:RelateYnnXnn} with eq. \eqref{eq:IncorrectOutputEq} used in widely used design tools, we see that there needs to be an additional correction needed in computing $y_{n|n}$ when process noise feeds into measurements, and traditional equality $y_{n|n} = C x_{n|n} + D u_{n}$ would not hold. This is so, because conditional expectation of process noise, knowing the measurement with its feed-through, is non-zero.
\begin{align} \label{eq:ExpectationWnGivenZn}
w_{n|n} \doteq \E(w_{n}|z_{n}) = K_{g2} \cdot (z_{n} - ( C_{m} x_{n|{n-1}} + D_{m} u_{n} ))
\end{align}
and eq. \eqref{eq:RelateYnnXnn} can be expressed as 
\begin{align}\label{eq:CorrecetdOutputEq}
y_{n|n} = C x_{n|n} + D u_{n} + H \E(w_{n|n})
\end{align}

Note that in absence of the feed-through and correlation between process and measurement noise ($H_{m}=0$ and $N=0$), $Kg_{2}$ would be zero and \eqref{eq:RelateYnnXnn} reverts to previous implementation \eqref{eq:IncorrectOutputEq}. Thus this correction only affects cases with correlated process and measurement noises or direct-feed through.

Note that the state-update and predicton equations from eq. \eqref{eq:XYXnextupdate} exactly matches with those from previous references and implementations. Thus the only modification proposed here is in computation of the output estimate, $y_{n|n}$.

\section{Steady-state Kalman filter}\label{sec:ssKF}
Above time-varying update equations simplify considerably for a steady-state Kalman filter.
The variance update equation \eqref{eq:VarupdateXNext} simplifies to the Riccati equation.
{\small
\begin{align}\label{eq:RiccatiDiscrete}
\begin{split}
&P = (APA^{T} + GQG^{T}) - \left(A P C_{m}^{T} + G Q H_{m}^{T} + G N \right) \cdot \\
& \hskip 2em \left(C_{m} P C_{m}^{T} + \bar{R} \right)^{-1} 
\left(A P C_{m}^{T} + G Q H_{m}^{T} + G N \right)^{T}
\end{split}
\end{align}
}
and the update equations \eqref{eq:XYXnextupdate} become linear time-invariant once $P$ is fixed, as $K_{g}$, $K_{g2}$ and all the coeffients become constant, and estimates evolve as per following dynamical system.
{\small
\begin{align} \label{eq:SteadyKF} 
\begin{split}
 x_{{n+1}|n} &= (A - M_{A,G} \,C_{m}) x_{n|{n-1}} +
\begin{bmatrix} M_{A,G} \\ B - M_{A,G} \,D_{m}\end{bmatrix} \begin{bmatrix} z_{n}\\ u_{n} \end{bmatrix} \\
\begin{bmatrix} x_{n|n}\\ y_{n|n}\end{bmatrix} &= 
\begin{bmatrix}
     I - K_{g} \,C_{m} \\ C - M_{C,H} \,C_{m}
\end{bmatrix} \cdot x_{n|{n-1}} \\
& \hskip 4em +
\begin{bmatrix} K_{g} && -K_{g}D_{m}\\ M_{C,H} && D - M_{C,H} \,D_{m}\end{bmatrix} \cdot 
\begin{bmatrix} z_{n} \\ u_{n}\end{bmatrix}
\end{split}
\end{align}
}

\section{Remarks on continuous time estimation}\label{sec:contKF}
As seen above, for the equations for discrete time Kalman filter, the only correction needed is in the formula for $y_{n|n}$ in eqs. \eqref{eq:RelateYnnXnn} and \eqref{eq:XYXnextupdate} to compute minimum variance output estimate.
The same correction carries over to the continuous time system. In continuous time systems, the state update is continuous, while measurements may be take at variable times. Thus prior state estimate $x_{n|n-1}$ is discrete time systems is replaced by prior estimate $\hat{x}^{-}_{t}$ and it's variance would be $\hat{P}^{-}_{t}$, and we can use the second equation from eq. \eqref{eq:XYXnextupdate} to derive continuous-time analogue of optimal posterior output estimate, $\hat{y}^{+}_{t}$.

\begin{align}   \label{eq:RelateYnnXnnCont}
\begin{split}
\hat{y}^{+}_{t} &= C \hat{x}^{-}_{t} + D u_{t}+ (C K_{g} + H K_{g2}) \\
&\hskip 8em \cdot \left(z_{t} - ( C_{m} \hat{x}^{-}_{t} + D_{m} u_{t} ) \right)
\end{split}
\end{align}
where $K_{g}=P^{-}_{t}C_{m}' (C_{m}P^{-}_{t}C_{m}' + \bar{R})^{-1}$, $K_{g2} = (Q H_{m}'+N) (C_{m}P^{-}_{t}C_{m}' + \bar{R})^{-1}$, which is same as the discrete time analog, except that $P^{-}_{t}$ is the prior variance of the prior state estimate $\hat{x}^{-}_{t}$.

\section{A Numerical Example}\label{sec:Example}
Consider a simple example of a system with a direct feed-through.
{\small
\begin{align*}
\dot{x} &= -0.1 x + 2 w, \hskip 1em
y = \begin{bmatrix}1\\0\end{bmatrix}x + \begin{bmatrix}1\\1\end{bmatrix}w,
\hskip 1em 
z = y(1) + v\\
Q &= \mathtt{var}(w) = 1, \hskip 1em 
R = \mathtt{var}(v) = 0.1, \hskip 1em
N = \mathtt{cov}(w,v) = 0.3
\end{align*}
}
Note that the first output $y(1)$ has feed-through from process noise $w$ and is measured with certain noise as $z$ and the second output $y(2)=w$ is an estimate of unknown input expressed as an auxiliary output.

Kalman estimators of three kinds were implemented on discretized version of above problem with time-step $0.1$ seconds. These were:
\begin{enumerate}[(a)]
\item Time varying Kalman filter as per eq. \eqref{eq:XYXnextupdate}, denoted by legend {\tt new}.
\item Steady state Kalman filter as per eq. \eqref{eq:SteadyKF}, denoted by legend {\tt new ss}.
\item Steady state Kalman filter using {\tt kalman.m} function in Matlab \cite{MatlabKalmanFnc2015b} controls toolbox, denoted by legend {\tt prev ss}.
\end{enumerate}

\subsection{Nominal performance}
\begin{figure}[h!]
  \includegraphics[width=\linewidth]{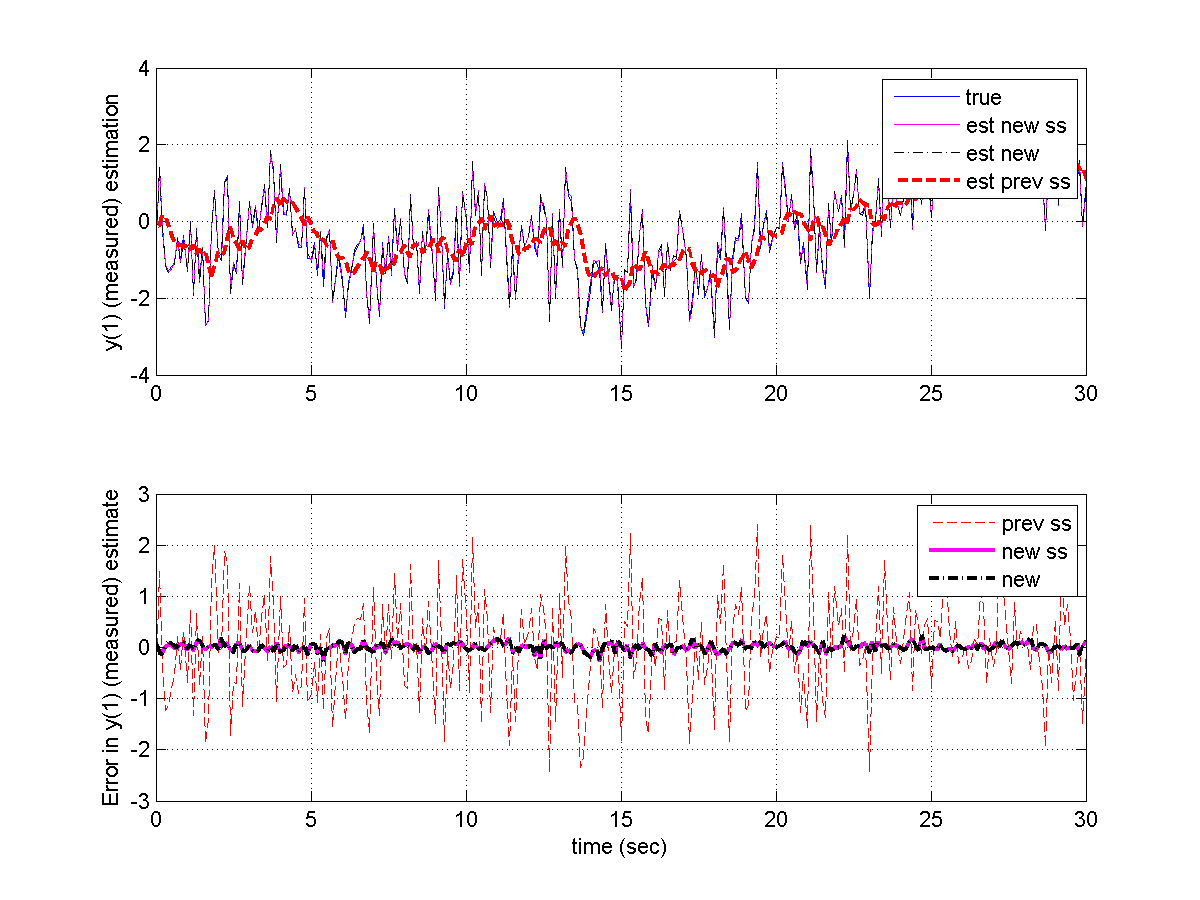}
  \caption{Measured output estimation $y(1)$}
  \label{fig:estErrY1}
\end{figure}

The estimates of measured output $y(1)$ by three methods are compared in fig. \ref{fig:estErrY1}, which also compares the estimate errors. It is clear that method (c) fails to capture the high frequencies present in the measured output, and only captures contribution due to state, which evolves smoothly. The estimate errors variances using both methods (a) and (b) is 0.0910, which matches eq. \eqref{eq:VarupdateYmeas}. This is 10 times lower  than method (c)  whose error variance is 0.99, which is close to $H_{m}QH_{m}^{T}+R$ as claimed in \cite{kwakernaak1972linear}.  The estimates from methods (a) and (b) are nearly identical after the initial transient. This matches with the common intuition that time-varying Kalman filter rapidly converges to steady-state Kalman filter, as the value of state variance $P_{n}$ approaches Riccati solution \eqref{eq:RiccatiDiscrete}.

\begin{figure}[h!]
  \includegraphics[width=\linewidth]{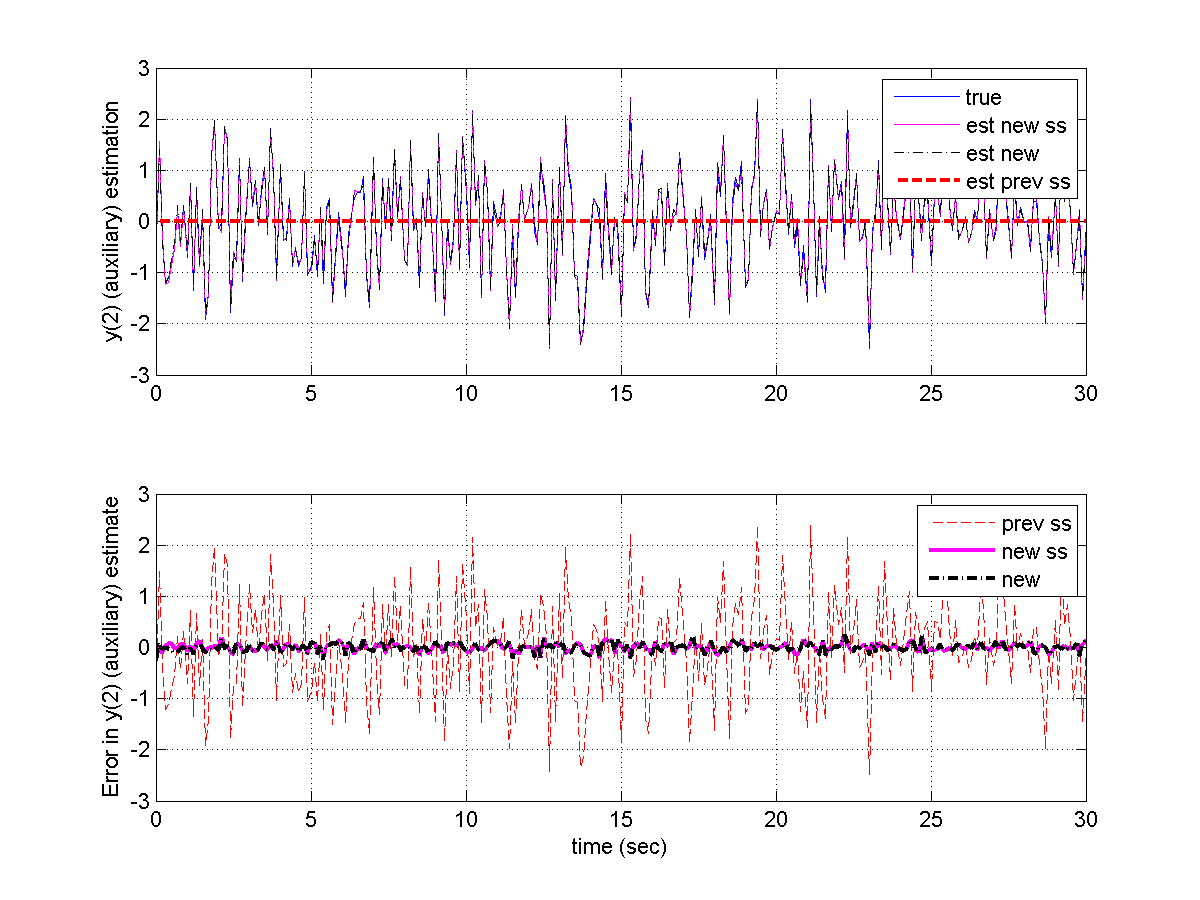}
  \caption{Unknown input estimation $y(2)$}
  \label{fig:estErrY2}
\end{figure}

 Fig. \ref{fig:estErrY2} compared the unknown input estimates by three methods, and again shows much better performance by the new method in unknown input estimation.

\begin{figure}[h!]
  \includegraphics[width=\linewidth]{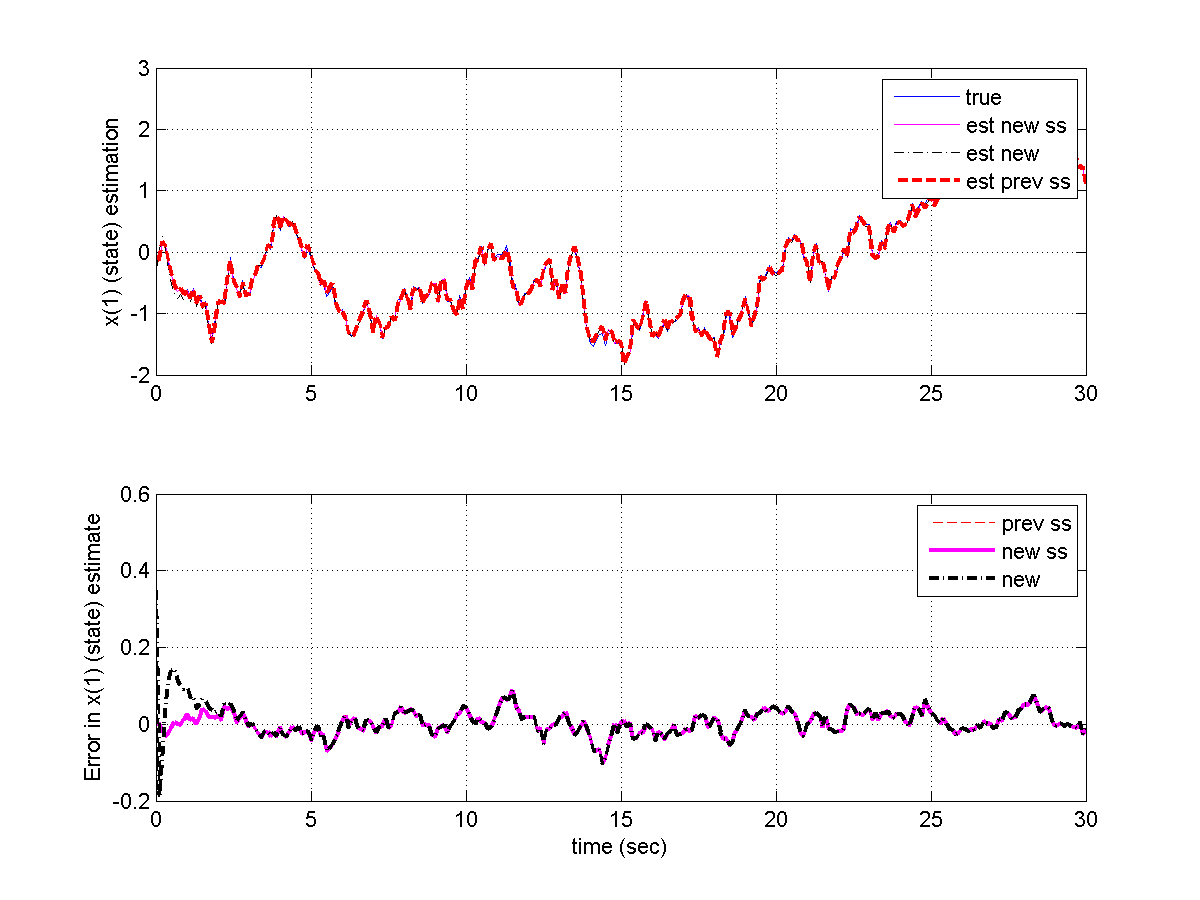}
  \caption{State estimation $x(1)$}
  \label{fig:estErrX1}
\end{figure}

 Fig. \ref{fig:estErrX1} compares the state estimates by three methods, and they show near-perfect match, which shows that state estimation by new method per eq. \eqref{eq:XupdateNext} is done in exactly same manner as method (c).

\subsection{Robustness to unmodeled bias in disturbance}
\begin{figure}[h!]
  \includegraphics[width=\linewidth]{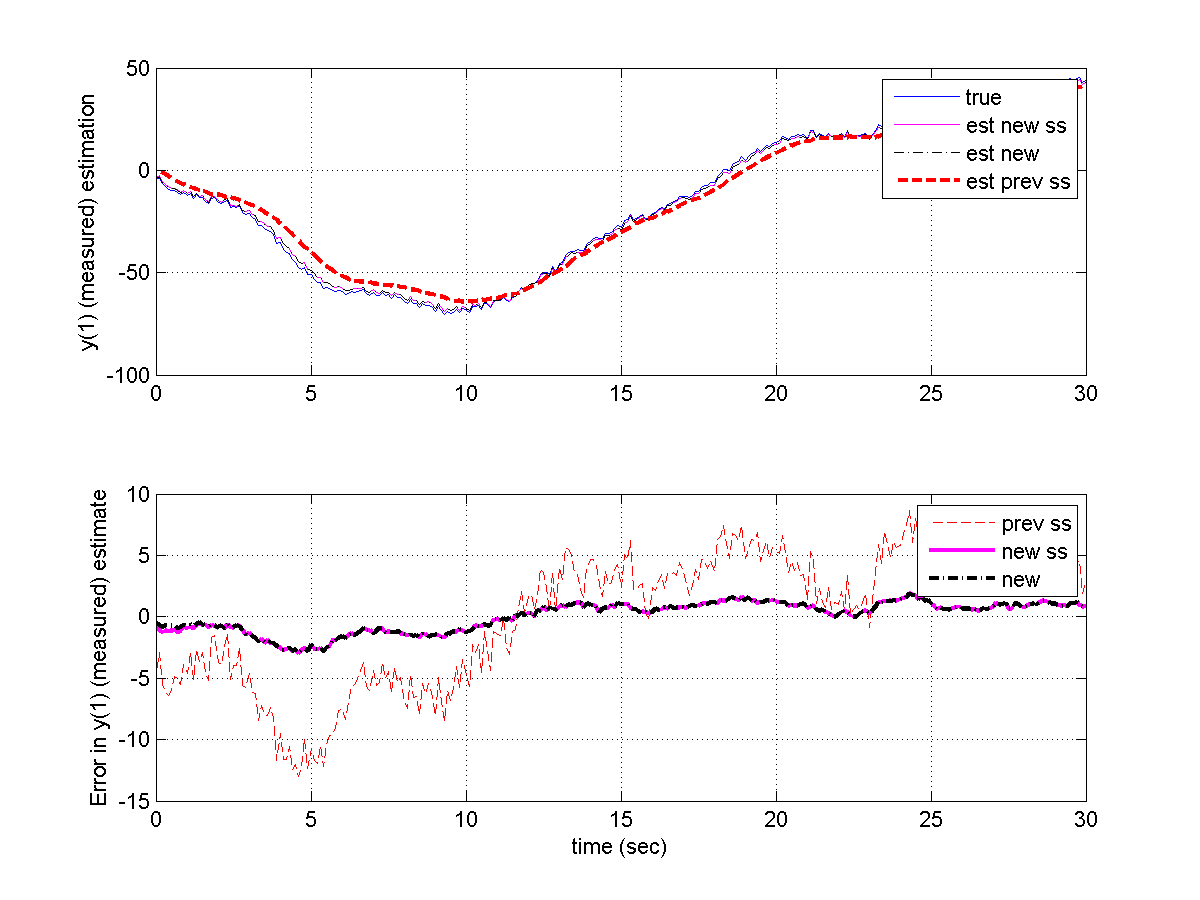}
  \caption{Measured output $y(1)$ estimation with unmodeled bias in process noise}
  \label{fig:estErrWithWBiasY1}
\end{figure}

\begin{figure}[h!]
  \includegraphics[width=\linewidth]{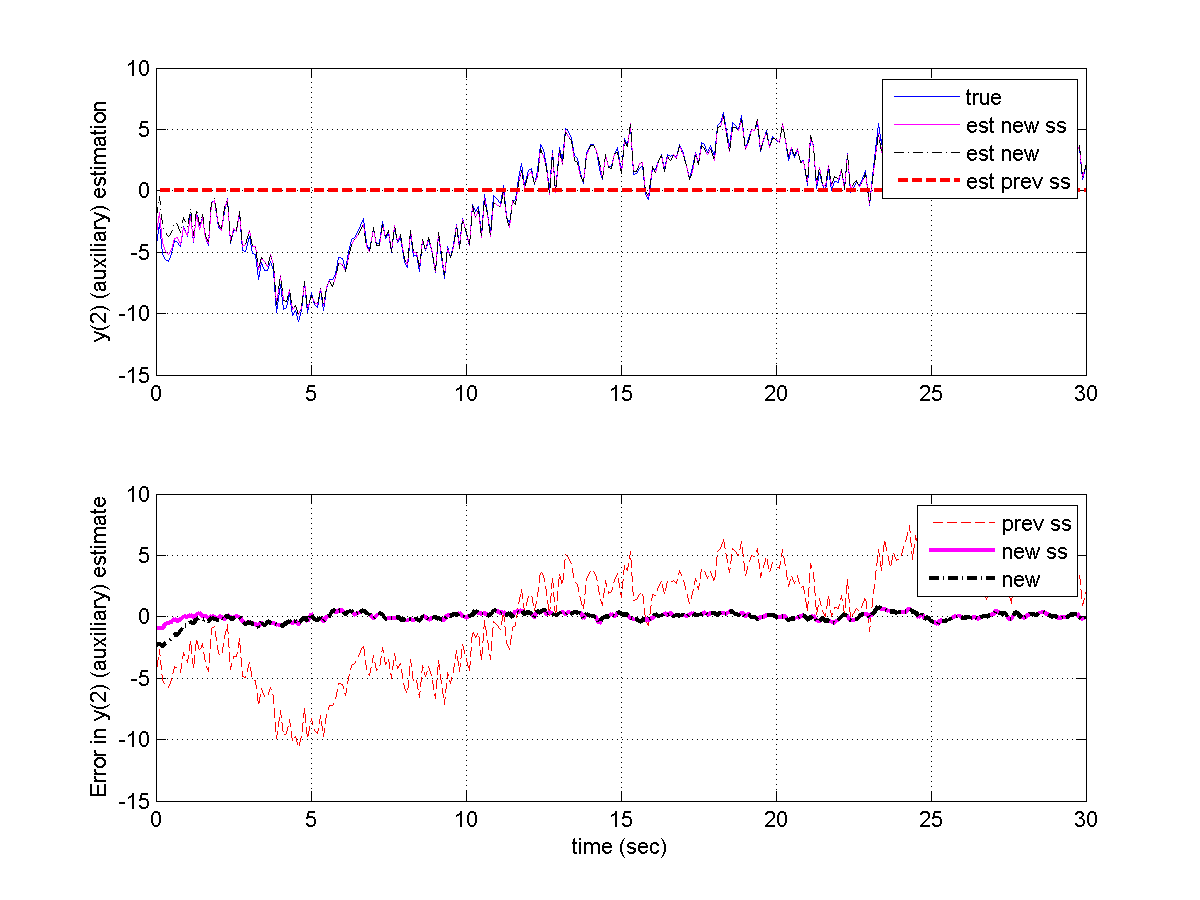}
  \caption{Unknown input $y(2)$ estimation with unmodeled bias in process noise}
  \label{fig:estErrWithWBiasY2}
\end{figure}

Usually it is hard to guarantee robustness of optimal estimator to unmodeled disturbance dynamics. Still, as an example, the three estimators above were tested against an unmodeled random walk drift in addition to the white noise in the disturbance $w$.  The {\tt true} signal in fig. \ref{fig:estErrWithWBiasY2} shows this disturbance. Fig. \ref{fig:estErrWithWBiasY1} compares $y(1)$ estimation and shows great improvement by new estimators over previous method (c), which has a clear drift in estimate error. Fig. \ref{fig:estErrWithWBiasY2} compares $y(2)$ or unknown input estimation. Since $y(2)=0\cdot x + w$, the estimate from method (c) is zero and estimate error is large. On the same problem, the new method estimates drift and high frequency disturbance quite well.


\section{Concluding Remarks}
In summary, until now, due to a gap in output estimate update rules in Kalman filter and sub-optimal implementation within widely used design tools (\cite{MatlabKalmanFnc2015b}, \cite{LabviewKalmanFnc}, \cite{MathematicaKalmanFnc}, \cite{MapleKalmanFnc}, \cite{OctaveKalmanFnc}) per eq. \eqref{eq:IncorrectOutputEq}, designers had no way of accommodating direct feed-through of even a part of the process noise into measurements, except by expanding state-space to include noise states, which may or may not be physical. But this workaround left much to be desired as it was constrained by smoothness requirements on the unknown input estimates and thus incurred large estimation errors at high frequencies and reduced the estimation bandwidth for such outputs. Bandwidth could only be increased at the cost of simplicity or robustness due to increased model order or new tuning parameters.

The proposed simple correction to output estimates per eq. \eqref{eq:RelateYnnXnn}, fills this gap as it makes use of posterior estimate of the unknown input computed from measurement which contain its feed-through. Thus it guarantees minimum variance output estimate. The performance of estimator has been demonstrated through a simple numerical example. The robustness to unmodeled inputs (e.g. bias or drift) seems to have greatly improved by this correction.

Such estimator is directly useful for the problems of estimating thrust from accelerometers mounted on flexible structures, especially if the thrust is broad-band and hard to model physically. The paper also makes the case for correcting the implementations in commonly used toolboxes for control system design, as all the surveyed implementations/documentations suffer from the sub-optimality in the output estimation for problems with direct feed-through or correlated process and measurement noises. 



\section*{ACKNOWLEDGMENT}
Many thanks to Cecilia Mazzaro, Xiongzhe Huang, Hullas Sehgal and Kirk Mathews from GE for their review which refined this work, and to Pascal Gahinet and Nirja Mehta from Mathworks, who reviewed and acknowledged and accepted this improvement and implemented it in Matlab release 2016a.
\bibliographystyle{IEEEtran}
\bibliography{IEEEabrv,rootbib}

\begin{thebibliography}{10}
\providecommand{\url}[1]{#1}
\csname url@rmstyle\endcsname
\providecommand{\newblock}{\relax}
\providecommand{\bibinfo}[2]{#2}
\providecommand\BIBentrySTDinterwordspacing{\spaceskip=0pt\relax}
\providecommand\BIBentryALTinterwordstretchfactor{4}
\providecommand\BIBentryALTinterwordspacing{\spaceskip=\fontdimen2\font plus
\BIBentryALTinterwordstretchfactor\fontdimen3\font minus
  \fontdimen4\font\relax}
\providecommand\BIBforeignlanguage[2]{{%
\expandafter\ifx\csname l@#1\endcsname\relax
\typeout{** WARNING: IEEEtran.bst: No hyphenation pattern has been}%
\typeout{** loaded for the language `#1'. Using the pattern for}%
\typeout{** the default language instead.}%
\else
\language=\csname l@#1\endcsname
\fi
#2}}

\bibitem{KalmanCounterExample}
\BIBentryALTinterwordspacing
Example script as a proof of suboptimality of traditional Kalman filter
  implementations. [Online]. Available:
  \url{http://arxiv.org/src/1510.00091v2/anc/testKFcorrection.m}
\BIBentrySTDinterwordspacing

\bibitem{MatlabKalmanFnc2016}
\BIBentryALTinterwordspacing
Documentation of Kalman filter design function in Matlab 2016 controls
  toolbox. [Online]. Available:
  \url{http://www.mathworks.com/help/control/ref/kalman.html}
\BIBentrySTDinterwordspacing

\bibitem{Kalman1960}
R.~E. Kalman, ``A new approach to linear filtering and prediction problems,''
  \emph{Transactions of the ASME--Journal of Basic Engineering}, vol.~82, no.
  Series D, pp. 35--45, 1960.

\bibitem{kwakernaak1972linear}
\BIBentryALTinterwordspacing
H.~Kwakernaak and R.~Sivan, \emph{Linear optimal control systems}, ser.
  Wiley-Interscience publication.\hskip 1em plus 0.5em minus 0.4em\relax Wiley
  Interscience, 1972. [Online]. Available:
  \url{https://books.google.com/books?id=mf0pAQAAMAAJ}
\BIBentrySTDinterwordspacing

\bibitem{lewis1986optimal}
F.~L. Lewis and F.~Lewis, \emph{Optimal estimation: with an introduction to
  stochastic control theory}.\hskip 1em plus 0.5em minus 0.4em\relax Wiley New
  York et al., 1986.

\bibitem{FranklinWorkmanPowell1997}
G.~F. Franklin, M.~L. Workman, and D.~Powell, \emph{Digital Control of Dynamic
  Systems}, 3rd~ed.\hskip 1em plus 0.5em minus 0.4em\relax Boston, MA, USA:
  Addison-Wesley Longman Publishing Co., Inc., 1997.

\bibitem{kailath2000}
T.~Kailath, A.~Sayed, and B.~Hassibi, \emph{Linear Estimation}, ser.
  Prentice-Hall information and system sciences series.\hskip 1em plus 0.5em
  minus 0.4em\relax Prentice Hall, 2000.

\bibitem{simon2006}
D.~Simon, \emph{Optimal State Estimation: Kalman, H Infinity, and Nonlinear
  Approaches}.\hskip 1em plus 0.5em minus 0.4em\relax Wiley-Interscience, 2006.

\bibitem{MatlabKalmanFnc2015b}
\BIBentryALTinterwordspacing
Documentation of Kalman filter design function in Matlab 2015b controls
  toolbox. [Online]. Available:
  \url{http://www.mathworks.com/help/releases/R2015b/control/ref/kalman.html}
\BIBentrySTDinterwordspacing

\bibitem{LabviewKalmanFnc}
\BIBentryALTinterwordspacing
Documentation of discrete Kalman filter function in labview. [Online].
  Available:
  \url{http://zone.ni.com/reference/en-XX/help/371894G-01/lvsim/sim\_disckalmanfilter/}
\BIBentrySTDinterwordspacing

\bibitem{MathematicaKalmanFnc}
\BIBentryALTinterwordspacing
Documentation of discrete Kalman filter function in Mathematica. [Online].
  Available:
  \url{https://reference.wolfram.com/language/ref/KalmanEstimator.html}
\BIBentrySTDinterwordspacing

\bibitem{MapleKalmanFnc}
\BIBentryALTinterwordspacing
Documentation of discrete Kalman filter function in Maple. [Online].
  Available:
  \url{http://www.maplesoft.com/support/help/MapleSim/view.aspx?path=ControlDesign/Kalman}
\BIBentrySTDinterwordspacing

\bibitem{OctaveKalmanFnc}
\BIBentryALTinterwordspacing
Documentation of Kalman filter design function in Octave-forge controls
  toolbox. [Online]. Available:
  \url{http://octave.sourceforge.net/control/function/kalman.html}
\BIBentrySTDinterwordspacing

\bibitem{CondBivarDistrib}
\BIBentryALTinterwordspacing
Conditional distribution of a multi-variate normal distribution. [Online].
  Available:
  \url{https://en.wikipedia.org/wiki/Multivariate\_normal\_distribution\#Conditional\_distributions}
\BIBentrySTDinterwordspacing

\end{thebibliography}

%
%

\end{document}